\documentclass[final]{siamltex}
\usepackage{graphics,graphicx,verbatim,upquote}
\newlength\myverbindent 
\makeatletter
\def\verbatim@processline{%
 \hspace{\myverbindent}\the\verbatim@line\par}
\makeatother
\def\C{{\bf C}}

\def\Im{\hbox{\rm\kern .3pt Im\kern 1pt}}
\def\Re{\hbox{\rm\kern .3pt Re\kern 1pt}}
\def\arg{\hbox{\rm\kern .3pt arg\kern 1pt}}
\def\Res{\hbox{\scriptsize\rm\kern .3pt Re\kern 1pt}}
\def\inv{[\kern .5pt 0,1]}
\def\minv{[-\infty,0\kern .5pt]}
\def\mminv{(-\infty,0\kern .5pt]}
\def\om{\omega}
\def\ft{\tilde f}
\def\omb{\bar\omega}
\def\Itail{I_{\scriptsize\rm tail}}
\def\Imiddle{I_{\scriptsize\rm middle}}
\def\Ihead{I_{\scriptsize\rm head}}
\def\bigint{[10^{-80},1]}\def\smallint{[10^{-20},1]}
\def\pput(#1,#2)#3{\noindent\smash{\raise#2pt\hbox to 0pt
   {\kern #1pt #3\hss}}\ignorespaces}

\title{Reciprocal-log approximation and planar PDE solvers}

\author{Yuji Nakatsukasa
and Lloyd N.~Trefethen\thanks{
\texttt{nakatsukasa@maths.ox.ac.uk} and \texttt{trefethen@maths.ox.ac.uk},
Mathematical Institute, University of Oxford, Oxford, OX2 6GG, UK.}}

\begin{document}

\maketitle


\begin{abstract}
This article is about both approximation theory and the numerical
solution of partial differential equations (PDE\kern .5pt s).  First we
introduce the notion of {\em reciprocal-log\/} or {\em log-lightning approximation}
of analytic functions with branch point singularities at points
$\{z_k\}$ by functions
of the form $g(z) = \sum_k c_k /(\log(z-z_k) - s_k)$,
which have $N$ poles potentially distributed along a Riemann surface.
We prove that the errors of best reciprocal-log approximations decrease
exponentially with respect to $N$ and
that exponential or near-exponential convergence
(i.e., at a rate $O(\exp(-C N / \log N))$)
also holds for near-best approximations with
preassigned singularities constructed
by linear least-squares fitting on the boundary.  We then apply
these results to derive a ``log-lightning method'' for numerical
solution of Laplace and related PDE\kern .5pt s in two-dimensional
domains with corner singularities.  The convergence is near-exponential,
in contrast to the root-exponential convergence for the original
lightning methods based on rational functions.
\end{abstract}


\pagestyle{myheadings}
\thispagestyle{plain}
\markboth{\sc Nakatsukasa and Trefethen}
{\sc Reciprocal-log approximation and planar PDE solvers}

\section{\label{introd}Introduction}

In this paper we introduce a new problem in approximation theory.
On the interval $\minv$, it is known
that functions $F(s)$ like $e^{a s}$ ($\Re a > 0$) with essential singularities at $s=-\infty$
can be approximated by rational functions $r(s)$ with exponential
convergence as a function of $n$,
the degree of the rational function.
Our starting point in section~\ref{sec2} is the change of variables $z = e^s$
from $s\in \minv$ to $z \in\inv$.
This gives us exponentially convergent approximations on $\inv$
of functions $f(z)$ like $z^a$ ($a>0$) with branch point singularities
at $z=0$ by functions of the form $g(z) = r(\log(z))$.  In the
generic case where the poles $\{s_k\}$ of
$r$ are distinct, the approximation can be written
in the partial fractions form
\begin{equation}
g(z) = c_0^{} + \sum_{k = 1}^n {c_k \over \log(z) - s_k}.
\label{form}
\end{equation}
We show that approximations of this kind can be computed by
linear least-squares fitting, and we note the equioscillatory
characterization of the error in the case of minimax (best
$\infty$-norm) approximations.  With both best and near-best
approximations, one encounters the startling property that the
oscillations typically cluster double-exponentially near the
singularity, readily coming as close (in theory) as a distance
of, say, $10^{-1000}$.   This is in contrast to the well-known
case of rational approximation, where the clustering is just
exponential~\cite{clustering}.

Section~\ref{sec3} generalizes the discussion to domains
in the complex \hbox{$z$-plane} with $m\ge 1$ singularities, typically at
corners.  Here we consider approximations of the generic form
\begin{equation}
g(z) = \sum_{j = 1}^m \kern 2pt \sum_{k=1}^{n_j}\kern 1pt
{c_{jk} \over \log(z-z_j) - s_{jk}} + p_0^{}(z),
\label{formsm}
\end{equation}
where $p_0^{}$ is a polynomial of degree $n_0^{}$ (possibly~0), with
total number of degrees of freedom $N = 1 + \sum n_j$.  The
least-squares method extends to this case too and again gives
exponentially or near-exponentially convergent approximations.
(By near-exponential, we mean at a rate $O(\exp(-C N /\log
N))$ with $C>0$.)  This and other claims are first explained
heuristically and illustrated numerically.  Then they
are established rigorously by a succession of theorems in
section~\ref{sec4}, deriving accuracy estimates from Hermite
integrals, and section~\ref{sec5}, showing how a problem with
several singularities can be decomposed into single-singularity
problems by Cauchy integrals over open arcs.

As is typical with computations based on redundant expansions,
numerical instability resulting from the use of ill-conditioned
matrices may be a problem; we stabilize our computations by the
Arnoldi orthogonalization technique presented in~\cite{VA}.
Section~\ref{secsing} considers additional 
fine points associated with the
``$10^{-1000}\kern 1pt$'' effect mentioned above, proposing
in particular that it may be a good idea to constrain the singular
part of a reciprocal-log approximation to approach zero more
rapidly at each singularity.

Up to this point, the paper is all approximation theory and
algorithms.  The motivation, however, is numerical solution of
PDE\kern .5pt s by a ``log-lightning'' analogue of the lightning
solvers for Laplace, Helmholtz, and biharmonic equations presented
in~\cite{stokes}, \cite{lightning}, and \cite{helm}.  A proof of
concept is briefly presented in section~\ref{sec7} to establish
that, as intended, these approximations can be used to derive
exponentially convergent new methods for solving PDE\kern .5pt
s in planar domains with corner singularities.

In a word, this is a paper about approximation by functions with
logarithmic branch points rather than just pole singularities.
What starts out as a seemingly arbitrary rearrangement of familiar
rational approximation theory turns out to have potentially
significant consequences for scientific computing.

A publication of Driscoll and Fornberg in 2001 proposed the
introduction of a logarithmic term to approximate functions with jumps, with
coefficients computed in the manner of
quadratic Hermite--Pad\'e approximation~\cite{df}.
Their technique is different from what is being proposed here, but
the use of logarithms to resolve singularities gives a family resemblance.

\section{\label{sec2}Approximation on \boldmath $\inv$}
A classic problem in approximation theory is the rational
minimax approximation of $F(s) = e^s$ on the interval
$\minv$.
To be precise, for each\/ $n$, let us consider rational functions $r(s)$
of degree $n$, meaning that $r(s)$ can be written as $p(s)/q(s)$
for degree $n$ polynomials $p$ and $q$.
In the generic case where
the poles are distinct and finite, $r$ can be written in partial fractions
form as
\begin{equation}
r(s) = c_0^{} + \sum_{k = 1}^n {c_k \over s - s_k}
\label{partfrac}
\end{equation}
with coefficients $\{c_k\}$ and poles $\{s_k\}$.
Of course, some rational functions have confluent poles (i.e., poles
of order $>1$) and
cannot be written in this form, and we shall sometimes make use
of such approximations.

Cody, Meinardus, and Varga showed in 1969~\cite{cmv} that
approximations $r(s)$ to $F(s) = e^s$ exist with exponential
convergence, meaning $\|F-r\| = O(C^{-n})$ as $n\to\infty$ for
some $C>1$, where $\|\cdot\|$ is the $\infty$-norm on $\minv$.
The optimal constant is $C = 1/H = 9.28903\dots,$ where $H$ is a
number related to elliptic integrals known as Halphen's constant;
for a discussion with references, see chapter~25 of~\cite{atap}.
Exponential convergence of best rational approximations at the same
asymptotic rate applies to $F(s) = e^{as}$ for any $a$ with $\Re
a > 0$ and to certain other functions $F(s)$ related to $e^{as}$
by rational transformations~\cite{ss}.  Moreover, convergence at
this rate occurs not just on $\minv$ but throughout the complex
plane~\cite{ss}.

By the change of variables $z = e^s$ and the definitions $f(z) = F(s)$ and
$g(z) = r(s)$, we transplant this theory to a new setting.
Now the domain is $z\in \inv$, and $f$ is a function such as
$z^a$ that may have a branch point singularity at $z=0$.
What is new is that (\ref{partfrac}) now takes the form
\begin{equation}
g(z) = c_0^{} + \sum_{k = 1}^n {c_k \over \log(z) - s_k},
\label{partfracr}
\end{equation}
where $\log$ denotes the standard branch of the logarithm.
Thus, by the transplantation of existing theory, we know
that a wide class of functions on $\inv$ with branch point
singularities at $z=0$ can be approximated by functions of the
form (\ref{partfracr}) with exponential convergence at the rate
$O((9.28903\dots)^{-n})$.

More important for scientific computing are approximations that
are not best but readily computable, still with exponential
convergence though not at the optimal rate.  One way to derive
such approximations is by expressing $F(s)$ as an integral along
a Hankel contour winding around $(-\infty,0\kern .5pt]$ in the
$s$-plane.  If the integral is approximated by a truncation of
the transplanted equispaced trapezoidal rule on the real axis, one
obtains rational approximations (\ref{partfrac}) to $F(s)$, hence
reciprocal-log approximations (\ref{partfracr}) to $f(z)$, which
converge at rates $O(C^{-n})$ with values of $C$ on the order of~3
or~4 depending on details of contours and parameters.  The leader
in developing these methods has been Weideman~\cite{tws,wt};
for a review, see section 15 of~\cite{trap}.

The approach we shall make use of is a variation on this theme
motivated by the
lightning PDE solvers of~\cite{stokes,lightning,helm} and
section~\ref{sec7}.\ \ Rather than deriving
coefficients from the trapezoidal rule,
we will obtain them by solving
a linear least-squares problem.  Typically this will be posed on
a discrete subset of $\inv$, with exponential
clustering of sample points at $z=0$.
Experiments show that an effective choice of singularities $\{s_k\}$ for
such approximations lies on a parabolic contour in the style of~\cite{tws,wt},
\begin{equation}
s_k = {n\over 4}(1+i \theta_k)^2, \quad \theta_k = -\pi + 2\pi
(k-\textstyle{1\over 2})/n, \quad 1 \le k \le n.
\label{hankelpts}
\end{equation}
For simplicity, however, we will also sometimes
use a configuration in which all the points $s_k$ are taken to be confluent at
a single value $s_0$ of scale $O(n)$; we normally take $s_0 = n/2$.
In this case (\ref{partfracr}) changes to the confluent form
\begin{equation}
g(z) = p\kern -1pt \left( {1\over \log(z)-s_0} \right),
\label{partfracrc}
\end{equation}
where $p$ is a polynomial of degree $n$.  As we shall prove in
section~\ref{sec4}, these methods give exponential convergence as
$n\to\infty$.

We illustrate these ideas with a MATLAB code segment whose fifth
line realizes the distribution (\ref{hankelpts}):

{\small
\begin{verbatim}

      M = 1000;
      z = logspace(-24,0,M)';
      f = sqrt(z);
      n = 10;
      kk = 1:n; tk = -pi + 2*pi*(kk-.5)/n; sk = (n/4)*(1+1i*tk).^2;
      A = [z.^0 1./(log(z)-sk)];
      c = A\f;
      g = A*c;

\end{verbatim}
\par}

\begin{figure}
\begin{center}
\vspace{4pt}
\includegraphics[scale=.96]{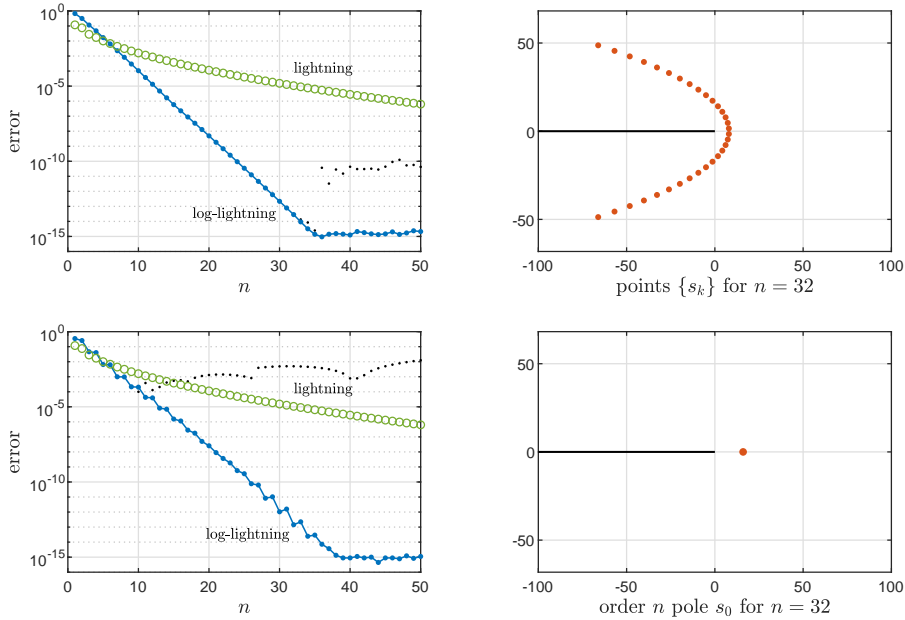}
\vspace{-2pt}
\end{center}
\caption{\label{fig1}Exponential convergence of approximations
$(\ref{partfracr})$--$(\ref{hankelpts})$ (top row)
and $(\ref{partfracrc})$ (bottom row) to $\sqrt z$ computed by least-squares
fitting on an exponentially clustered grid 
in $\inv$.  (These data are obtained with
Arnoldi stabilization~{\rm\cite{VA}}; the numbers obtained without
stabilization are shown by the small dots.)
A strong contrast is apparent with the root-exponential
convergence of lightning rational approximations~{\rm \cite{lightning}}
for the same problem.
The right column shows
the points $\{s_k\}$ or $s_0$ in the $s$-plane for $n=32$.  Since the imaginary
parts $\Im s_k$ extend far beyond $[-\pi,\pi]$ in
the first case, the values\/ $\exp(s_k)$, which are the poles of $g(z)$, lie on
many different sheets of a Riemann surface.}
\end{figure}

\noindent 
The maximum error on $\inv$ is $1.5\times 10^{-4}$, and this number
decreases exponentially as $n\to\infty$, as shown in the first
image of Figure~\ref{fig1}.  There is a marked contrast with the
root-exponential behavior of the lightning rational approximation
of~\cite{lightning}.  The second image shows the singularities
$\{s_k\}$ in the case $n=32$.  Note that although the geometry
in the $s$-plane is just an arc, the corresponding points $z_k
= \exp(s_k)$ in the $z$ domain lie along a spiral on the Riemann
surface of $g(z)$ associated with its definition in terms of
the function $\log(z)$,
clustering exponentially at both ends at $z=0$.  With
$n=32$, the last value in the upper half $s$-plane is $s_{32}^{}
\approx -66.1 + 48.7i$.  The corresponding value $z_{32}^{}$ has
modulus $1.0\times 10^{-31}$ and argument $15.5\pi$, putting it
8 sheets above the main branch of the Riemann surface of $g(z)$.

The second row of Figure~\ref{fig1} shows corresponding results
for the confluent-pole approximation (\ref{partfracrc}) with $s_0
= n/2$.  These poles are further from optimal, and the exponential
convergence is about 15\% slower.  Certain fine points associated
with confluent computations are discussed in section~\ref{secsing}.
To minimize the effects very near the singularity discussed there and
make the comparison of Figure~\ref{fig1} a fair one, the least-squares
fitting for this figure was done on a grid of $2000$ points logarithmically
spaced in $[10^{-100},1]$.

\begin{figure}
\begin{center}
\vspace{4pt}
\includegraphics[scale=.76]{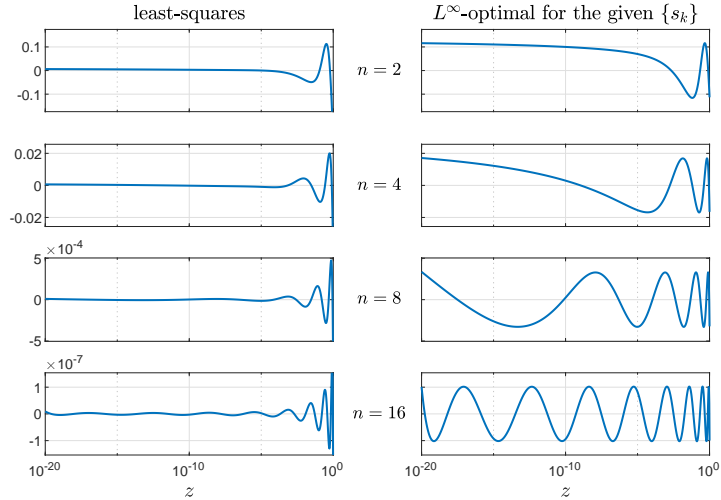}~~
\vspace{-4pt}
\end{center}
\caption{\label{fig2}Error curves for reciprocal-log approximations
$(\ref{partfracr})$
of $\sqrt z$ on $\inv$ with $n = 2, 4, 8, 16$  with singularities
$s_k$ given by $(\ref{hankelpts})$.  In the left column, coefficients
$\{c_k\}$ are computed
by least-squares fitting on an exponential clustered
grid from of $1000$ points from $10^{-20}$ to $1$.  In the
right column, the fits are then improved to $L^\infty$ optimal for the
given $\{s_k\}$ by a Lawson iteration (iteratively reweighted
least-squares)~{\rm \cite{lawson}}.  This improves the maximal error
by about a factor of\/~$2$ at the cost of eliminating the enhanced
accuracy for $z\approx 0$.
Each pair shares a common vertical scale, with
each curve showing $n+1$ extrema of alternating sign.}
\end{figure}

\begin{figure}
\begin{center}
\vspace{4pt}
\includegraphics[scale=.76]{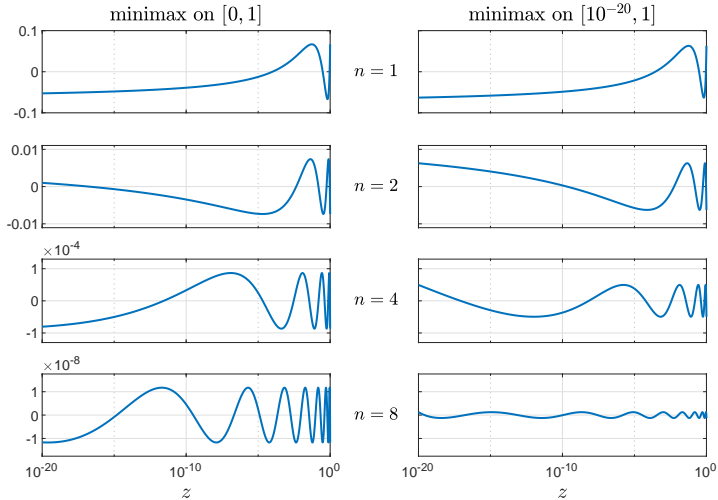}
\vspace{-4pt}
\end{center}
\caption{\label{fig3}Error curves for true minimax 
reciprocal-log approximations 
of $\sqrt z$ on $\inv$ with $n = 1, 2, 4, 8$ on $\inv$ (left)
and $[10^{-20},1]$ (right).  Note that these values of $n$ are
half those of Figure~$\ref{fig2}$, reflecting the roughly
squared accuracy of minimax approximations.
Again each pair shares a common vertical scale, and now there
are $2n+2$ equioscillating extrema.}
\end{figure}

Figure~\ref{fig2} shows error curves $g(z) - \sqrt z$, $z\in
[10^{-20},1]$, for approximations (\ref{partfracr}).  On the
left these are the approximations computed by least-squares
fitting with the singularities $\{s_k\}$ of (\ref{hankelpts}),
as in the first row of Figure~(\ref{fig1}).  Note that the
maximal error decreases exponentially with $n$ and that the
error is much smaller for $z\approx 0$ than for $z\approx 1$,
because the clustered grid puts more weight there.
The right column shows error curves for approximations with
the same $\{s_k\}$ whose coefficients have been
modified to achieve $L^\infty$ optimality by 20
steps of a linear Lawson iteration, also known as iteratively
reweighted least-squares~\cite{lawson}.
%
%
%
These equioscillatory error curves are elegant, but we doubt
whether the ``optimal'' approximations are superior in practice.
Their maximal error is lower by only about a factor of 2, at
the cost of greatly worsening the accuracy for $z\approx 0$.
See chapter~16 and Myth~3 of Appendix~A in~\cite{atap}.

The poles of our approximations $g(z)$ lie exponentially close
to the origin in the $z$-plane, at least when the formula
(\ref{hankelpts}) is used.  Exponential clustering of poles
is familiar in rational approximation theory~\cite{clustering}.
Another phenomenon appears in Figure~\ref{fig2} that is unfamiliar,
however, and that is {\em doubly-exponential} clustering of
oscillating extrema.  To see this, consider the leftmost minimum of
each error curve in the right column of the figure.  With $n=2,$
$4$, and $8$, these minima lie at approximately $z = 10^{-1.5}$,
$10^{-4.5}$, and (as shown by computations on a longer interval)
$10^{-25}$.  A rough model of such behavior would be $z_{\min{}}
= \exp(-1.7^{\kern .3pt n})$, which suggests values for $n=16$
and 32 on the order of $10^{-2000}$ and $10^{-10{,}000{,}000}$.
In the figure, we see nothing like this, since the computations
are carried out on $[10^{-20},1]$, and it would be impossible to
resolve much beyond $[10^{-300},1]$ in standard IEEE floating
point arithmetic.  However, this probably does not matter, as
the approximants on $\inv$ will not differ much from
those associated with a truncated interval such as $[10^{-20},1]$.
We discuss these matters further in section~\ref{secsing}.

All these approximants are linear, based on poles $s_k$
fixed a priori.  What about true minimax approximations, with
$\{s_k\}$ chosen optimally to minimize the error?  For the
current setting of $\inv$, everything is known about such
approximations since they are transplants of minimax rational
approximants on $\minv$~\cite[chap.~24]{atap}.  As illustrated
in Figure~\ref{fig3}, which was computed by the method outlined
on pp.~217--219 of~\cite{atap}, the errors will be about the
squares of what we have seen before, and now there are $2n+2$
equioscillatory extreme points.  We shall not discuss such
approximations further, because they become impractical once one moves
to planar domains.

\section{\label{sec3}Approximation on a planar domain}

We now turn to domains in the complex plane.  In the basic case
one has a closed Jordan region $E$ with piecewise analytic boundary $\Gamma$ and
a function $f$ continuous in $E$ and analytic in the interior.
We suppose that $f$ is analytic on $\Gamma$ except for branch point
singularities at a finite
set of points $\{z_j\}$, $1\le j \le m$, which may for
example be corners.  
As mentioned in the introduction,
we consider approximations to $f$ of the generic form
\begin{equation}
g(z) = \sum_{j = 1}^m \kern 2pt \sum_{k=1}^{n_j}\kern 1pt
{c_{jk} \over \log(z-z_j) - s_{jk}} + p_0^{}(z),
\label{formsagain}
\end{equation}
where $p_0^{}$ is a degree $n_0^{}$ polynomial,
with $N = 1 + \sum n_j$.
Each function \hbox{$\log(z-z_j)$} denotes a fixed choice
of a branch of the logarithm that is continuous in $E$.
Since standard hardware puts the branch of the logarithm on
the negative real axis, it is often convenient in a computer code to
work with rotated terms of the form $\log(e^{i\theta_j} (z-z_j))$ for
certain real numbers $\theta_j$.
We also employ the confluent variant
\begin{equation}
g(z) = \sum_{j = 1}^m \kern 2pt p_j\kern -2pt \left(
{1 \over \log(z-z_j) - s_j}\right) + p_0^{}(z),
\label{formsagainc}
\end{equation}
where each $p_j$ is a polynomial of degree $n_j$.
In both (\ref{formsagain}) and (\ref{formsagainc}), we often
take $n_0=0$, so that $p_0^{}$ reduces to a constant.

The essential point, to be proved in the next two sections, is that
exponential or near-exponential
convergence is guaranteed so long as the numbers
$n_j$ (for $j\ge 1$) increase in proportion to $N$ as $N\to\infty$, in the
absence of rounding errors of course.  In this section
we will explore just one example defined on the domain shown
in Figure~\ref{figphase}, a two-thirds bite of the
unit disk.  The function is
\begin{equation}
f(z) = z \log(-\textstyle{1\over 2}z)\cdot  (1-z/\om)^{1/2} \cdot
(1-z/\omb)^{3/2}
\label{3corner}
\end{equation}
with $\om = e^{\pi i /3}$, with branch point singularities at
each of the corners.  The boundary $\Gamma$ is discretized by
500 points along each segment exponentially clustered down to
distances from the corners of about $10^{-14}$.

\begin{figure}
\begin{center}
\vspace{7pt}
~\hbox{\kern .3in}\includegraphics[scale=.83]{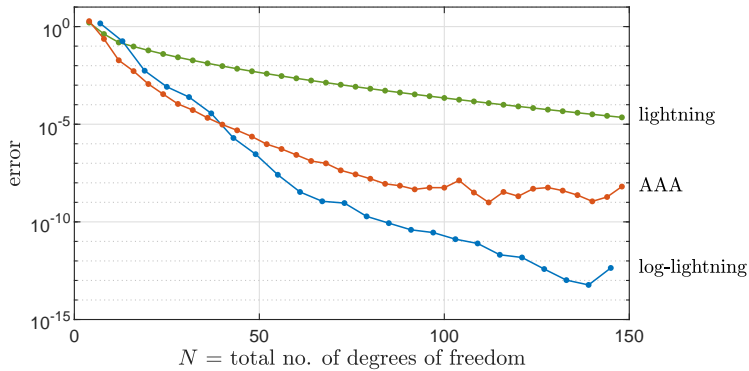}
\vspace{-4pt}
\end{center}
\caption{\label{figconv} Convergence of three approximants to
the function\/ $f$ of\/ $(\ref{3corner})$ on the domain $E$
shown in Figure~$\ref{figphase}$.  Down to accuracy of\/ $10^{-8}$ or so,
one sees root-exponential convergence of the AAA and lightning
approximations and exponential convergence of
the log-lightning approximations. Beyond this point, rounding errors
muddy the water but the log-lightning approximations are by far the
most accurate.}
\end{figure}

As indicated in Figure~\ref{figconv}, we computed
three approximations to $f$ for values of $N$ up to $150$.
One curve shows the accuracy of
AAA rational approximations of degree $N-1$.  (AAA approximations
are rational approximations that are computed very quickly by
an algorithm described in~\cite{AAA}, implemented in the
{\tt aaa} command of Chebfun; we calculated them 
with the ``cleanup'' option turned off.)  The AAA approximations are
close to optimal when the degree is not too large, and this curve shows
root-exponential convergence down to an error on the order of $10^{-8}$. 
Another curve shows lightning rational approximations~\cite{lightning} of the form
\begin{equation}
r(z) = \sum_{j = 1}^3 \kern 2pt \sum_{k=1}^{n_j}\kern 1pt
{c_{jk} \over z - s_{jk}} + p_0^{}(z)
\label{lightning}
\end{equation}
with $\sum n_j = N -1$, where $p_0^{}$ is a polynomial of
degree $n_0^{}$\kern .5pt; we take $n_0=N/4 - 1$ and $n_j = N/4$
for $j=1,2,3$.  
(The choice $n_0^{}=0$, which is often a good one for
(\ref{formsagain}) and (\ref{formsagainc}),
tends to slow down convergence significantly for (\ref{lightning}).)
The poles were fixed at distances $d_{jk} =
\exp(4(\sqrt{\vphantom{n_j}k} - \sqrt{\vphantom{k}n_j}\kern 1pt))$,
$1 \le k \le n_j$ from the corners $z_j$~\cite[eq.~(3.2)]{lightning}.
Note that the convergence is again root-exponential, at a rate five
or six times slower than with AAA with its adaptively determined
poles.  Clearly AAA has some advantages, but it has drawbacks of
sometimes greater sensitivity to rounding errors and occasional
placement of poles inside a domain where they are not wanted.
More fundamentally, as we shall discuss in section~\ref{sec7},
lightning methods generalize immediately to approximation of the
real part for solving Laplace problems, whereas AAA does not.

\begin{figure}
\vspace{1pt}
~~~~~~~~\includegraphics[scale=.83]{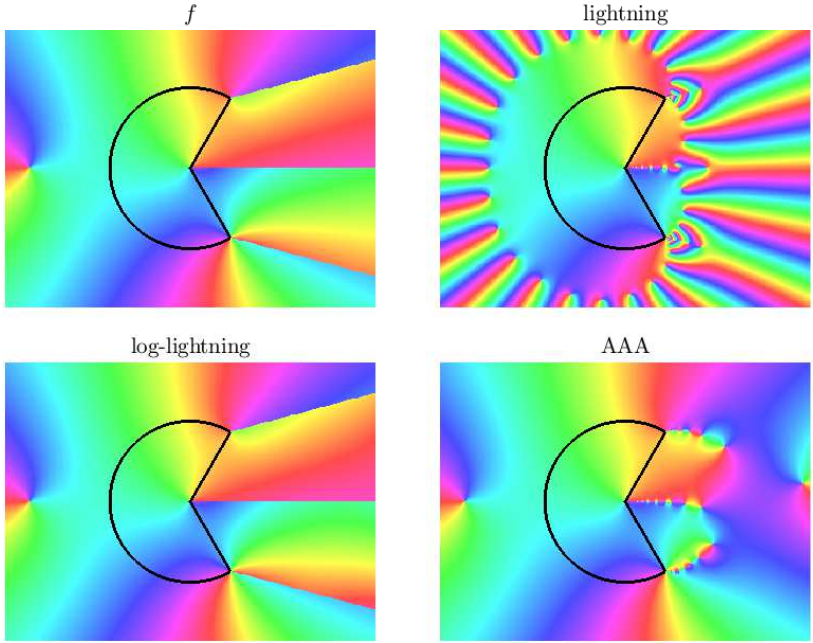}
\vspace{-30pt}
\caption{\label{figphase} Phase portraits of the function $f$
of\/ $(\ref{3corner})$ and three approximants with $N\approx 100$
degrees of freedom.  All
the approximations give many digits of accuracy in~$E$.
Outside~$E$, the branch cuts of\/ $f$ are approximated
by strings of poles for lightning and AAA approximants,
and by true branch cuts for the log-lightning approximant.
The circus tent effect in the lightning plot
reflects the polynomial term\/ $p_0^{}$ in\/ $(\ref{lightning})$.}
\vspace{-10pt}
\end{figure}

The final curve in Figure~\ref{figconv} shows log-lightning
approximations (\ref{formsagainc}) with $s_j = n_j/3$, now with
$n_0 = 0$ and $n_j = (N-1)/3$ for $j = 1,2,3$.
The convergence appears exponential
to around $10^{-9}$, after which it slows for reasons at least
partly related to rounding errors.

Figure~\ref{figphase} displays phase portraits~\cite{wegert}
for $f$ and its three approximations with $N=100$.  The rational
approximations simulate branch cuts by strings of poles, whereas
the reciprocal log approximation incorporates true branch cuts.
For $f$ itself, the shapes of the branch cuts are arbitrary.  If we
had coded the function differently, their directions could have
been altered, for example, but this would have had no effect on the
values of $f$ on $\Gamma$ and therefore on the three approximations
and the other three images.  For the lightning and log-lightning
approximations, the forms of the branch cuts are fixed by our
computation: they have been set to be straight lines bisecting
the exterior angles, the natural default choice in the absence
of other information.  For the AAA approximation, the branch cuts
are presumably more nearly optimal, curved in a fashion analyzed
by Stahl for the case of Pad\'e approximation~\cite{stahl}.

\section{\label{sec4}Convergence theorems: one singularity}

In this section we consider the reciprocal-log approximation
of a function $f(z)$ with a branch point singularity at $z=0$.
By arguments adapted from those of section~2 of~\cite{lightning},
we show that if $f$ can be analytically continued along any contour
in the complex plane avoiding $z=0$, then exponentially convergent
approximations exist with the singularities $\{s_k\}$ extending a
distance $O(n)$ into the right half-plane (Thm.~\ref{thm1}).  If
the analyticity assumption holds just in a bounded neighborhood of
$0$, then near-exponentially convergent approximations exist with
$\{s_k\}$ of bounded positive real part (Thm.~\ref{thmbounded}).

To make these ideas precise the theorems of this and the next
section use the following definition, which covers the familiar
singularities of type $z^a (\log z)^b$~\cite{lehman,wasow}.
The first condition, (\ref{taucond}), essentially asserts that
$f$ can be analytically continued along arbitrary contours in a punctured
neighborhood of each singularity $z_k$.  But all we really need
is continuation along curves that spiral in logarithmically to
$z_k$, so that after the $s = \log(z-z_k)$ change of variables
they lie in a $V$-shaped wedge in the $s$-plane, and that is why
(\ref{taucond}) takes the form it does.  The second condition,
(\ref{bpdef}), asserts that $f$ is H\"older continuous at $z_k$,
so that with $s=\log(z-z_k)$ we get exponential decay as $\Re s
\to -\infty$.  A function like $(z-z_k)^a \log(z-z_k)$, however,
will grow slowly as $z$ winds around $z_k$, and that is why
(\ref{bpdef}) includes the factor involving $\arg(z-z_k)$.

\medskip

\begin{definition}
\label{def}
Let $E$ be a closed set in the complex plane.  A
function $f$ is {\em analytic on $E$ with branch point
singularities at $z_1,\dots ,z_m$} on the boundary
if\/ $f$ is analytic in the interior of $E$ and
can be analytically continued to a neighborhood of
each boundary point of $E$ that is not among the $z_k$;
and if moreover there is a neighborhood of each $z_k$ within which $f$ can be
analytically continued to a multivalued function $\ft$ along 
any curve that avoids $z_k$ and satisfies
\begin{equation}
|\arg (z-z_k)| \le |z-z_k|^{-\tau}
\label{taucond}
\end{equation}
for some constant $\tau > 0$, with $\ft$ satisfying
\begin{equation}
|\ft(z) - \ft(z_k)| \le c \kern 1pt |z-z_k|^a (1+|\arg (z-z_k)|^b)
\label{bpdef}
\end{equation}
for some constants $a,b,c>0$.
\end{definition}
\medskip

For our first theorem, it is enough to use confluent singularities
and consider 
approximations of the form (\ref{partfracrc}) with $s_0= \sigma n$,
\begin{equation}
g(z) = p\kern -1pt \left( {1\over \log(z)-\sigma n} \right),
\label{e1}
\end{equation}
where $p$ is a polynomial of degree $n$ and $\sigma > 0$ is
a constant.  This result is sharp in the sense that faster than
exponential convergence is in general not possible, as we know from
the equivalence of the problem of reciprocal-log approximation
of $z^a$ on $\inv$ to that of rational approximation of $e^{as}$
on $\minv$.  The condition (\ref{h2}) looks almost like a repeat
of (\ref{bpdef}), but (\ref{bpdef}) is a local condition, only
needing to apply for $z$ near $z_k$, whereas (\ref{h2}) applies
even as $z\to\infty$.  With $s=\log z$, it is designed to ensure
that $F(s)$ is under control in the right half of the $s$-plane
as well as the left.

\medskip
\begin{theorem}
\label{thm1}
Let\/ $E$ be a compact set in 
$\C\kern 1pt \backslash\kern 1pt (-\infty,0\kern .5pt )$ with
$0$ on the boundary; more
generally, $E$ might be a multilevel surface wrapping around\/~$0$
a finite number of times.
Let\/ $f$ be an analytic function on $E$ with a branch point
singu\-larity at $z=0$ in the sense of
Definition~$\ref{def}$, and assume further that $f$ can be analytically
continued to a multivalued function $\ft$ along
all curves in the $z$-plane avoiding $z=0$, satisfying
in addition to $(\ref{bpdef})$
\begin{equation}
|\ft(z)| \le C \kern .3pt (1+|z|^{a'}) (1+|\arg z|^{b'})
\label{h2}
\end{equation}
for some constants $a',b',c' > 0$.
Then for all sufficiently small $\sigma>0$,
there exist functions of the form $(\ref{e1})$ satisfying
\begin{equation}
\|\kern .7pt f - g\| \le \exp(- C n ) \quad \forall n
\label{e3}
\end{equation}
for some\/ $C>0$, where $\|\cdot\|$ is the supremum norm on $E$.
\end{theorem}
\medskip

By the change of variables $s = \log z$ and subtraction of
$f(0)$, Theorem 1 is equivalent to the following theorem on
rational approximation.  As just mentioned, (\ref{h3}) constrains
$F(s)$ in the left half $s$-plane, and (\ref{h4}) in the right.
The function $e^s + s\kern .3pt e^{2s}$, for example, satisfies
the conditions with $a< 1$, $a'> 2$, and $b, b' > 1$.

\medskip
\begin{theorem}
\label{thm2}
Let\/ $E$ be a subset of\/ $\C$ with bounded\ $\Re E$
and\/ $|\Im E|$.
Let\/ $F$ be an entire function of\/ $s$ satisfying
for all sufficiently small $\Re s$
\begin{equation}
|F(s)| \le c \kern 1pt e^{a\kern .7pt \Res s} (1+|\Im s|^b)
\label{h3}
\end{equation}
and for all $s$
\begin{equation}
|F(s)| \le c' (1+ e^{a'\kern .7pt \Res s}) (1+|\Im s|^{b'})
\label{h4}
\end{equation}
for some constants $a, b, c, a', b', c'>0$.
Then for all sufficiently small\/ $\sigma>0$,
there exist type $(n-1,n)$ rational functions\/ $r$ of the form
\begin{equation}
r(s) = p\kern -1pt \left( {1\over s - \sigma n}\right),
\label{e4}
\end{equation}
where $p$ is a polynomial of degree $n$, with
\begin{equation}
\|F - r\| \le \exp(- C n ) \quad \forall n
\label{e5}
\end{equation}
\end{theorem}
for some\/ $C>0$, where $\|\cdot\|$ is the supremum norm on $E$.

\medskip

Thus our job is to prove Theorem~\ref{thm2}.  The proof
will be based on the Hermite integral formula for rational
interpolation~\cite[Thm.~2, Chap.~8]{walsh}.  Neither the
poles $\{s_k\}$ nor the interpolation points $\{\alpha_k\}$ in
the statement below need to be distinct, and an interpolation
point of multiplicity $\nu>1$ is interpreted, as usual, to mean
interpolation at that point of $f, f', \dots , f^{(\nu-1)}$.  By a
Hankel contour, we mean a continuous curve winding counterclockwise
around $(-\infty,0\kern .5pt]$ from $\infty$ to $\infty$.

\medskip
\begin{lemma}
\label{hif}
(Hermite integral formula)
Let $F$ be an analytic function of\/ $s$ satisfying\/ $(\ref{h3})$
on and inside a Hankel contour\/ $\Gamma$ in the $s$-plane.
Let interpolation points $\alpha_1,\dots, \alpha_n\in \minv$
and poles $s_1,\dots,s_n$ anywhere in
the $s$-plane be given; neither $\{\alpha_j\}$ nor $\{s_j\}$ need to be
distinct.
Let\/ $r$ be the unique type $(n-1,n)$ rational function with
poles at $\{s_j\}$
that interpolates $F$ at $\{\alpha_j\}$.  Then for any\/ $s$ enclosed
by\/~$\Gamma$,
\begin{equation}
F(s) - r(s) = {1\over 2\pi i} \int_\Gamma {\phi(s)\over \phi(t)}
{F(t)\over t-s} \kern 1pt dt,
\label{herm}
\end{equation}
where
\begin{equation}
\phi(s) = \prod_{j=1}^n  {s-\alpha_j\over s-s_j}.
\label{e7}
\end{equation}
\end{lemma}

\begin{figure}
\begin{center}
\vspace{10pt}
\includegraphics[scale=.62]{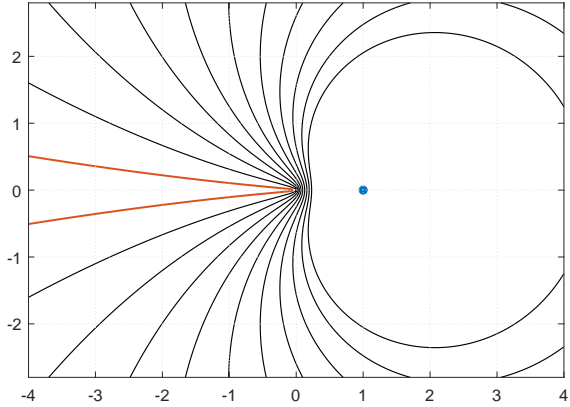}
\end{center}
\caption{\label{figcont}Level curves $n^{-1} \log |\phi(n\kern .5pt s)|
= 0,0.1,\dots, 1$
for $n=100$ to illustrate Lemma~$\ref{philemma}$.
As\/ $n$ increases, the contours converge to radial lines,
with the 
level $0$ contour (red) narrowing toward $\mminv$.
Thus $|\phi(n\kern .5pt s)|^{1/n}$ approaches
a value ${>}\kern 1pt 1$ at each point $s\in
\C \kern 1pt \backslash \{ (-\infty,0\kern .5pt] \cup \sigma \}$
and is bounded by $1$ for $s \in \mminv$.}
\vspace{1.5cm}
\end{figure}

In proving Theorem~\ref{thm2}, we will suppose
$E = \mminv$ without loss of generality;
the estimates change only by constant factors for other domains
since the arguments involve a geometry that scales with
$n$ as $n\to\infty$.
We will use this choice of interpolation points in $\mminv$:
\begin{equation}
\alpha_j = n\kern 1pt \sigma \kern -1pt\left({1+t_j\over 1-t_j}\right)^2, \quad 
t_j = e^{i\pi j/(n+1)}, \quad 1\le j \le n.
\label{interppts}
\end{equation}
These are derived from potential theory~\cite{levsaff,walsh}, where
poles and interpolation points are interpreted as point charges
of opposite signs, and a minimal-energy charge configuration
gives asymptotically optimal approximations as $n\to\infty$.
The function
\begin{equation}
s = n \kern .7pt \sigma \kern -1pt\left({1+t\over 1-t}\right)^2
\label{potent}
\end{equation}
maps the unit disk in the $t$-plane to the $s$-plane slit along
$\mminv$, with $t=0$ mapping to $s=n\kern .5pt \sigma$.  Thus the
equilibrium distribution of interpolation points on the unit circle
corresponding to poles at $t=0$, namely the uniform distribution,
maps to the distribution determined by (\ref{potent}) in the
$s$-plane slit along $(-\infty,0\kern .5pt]$ with poles at $s=n
\kern .5pt \sigma$, as in the function (\ref{e4}).  (The $\{\alpha
_j\}$ are called {\em Fej\'er--Walsh points}~\cite{starke}.)
The following lemma, illustrated in Figure~\ref{figcont}, was
provided to us by Peter Baddoo.

\medskip
\begin{lemma}
\label{philemma}
With $\alpha_j$ defined by $(\ref{interppts})$, 
$s_j = n\kern .7pt \sigma$, and\/ $s$ related to\/ $t$ by
$(\ref{potent})$, the function
$\phi(s)$ of $(\ref{e7})$ can be written
\begin{equation}
\phi(s) = {t^{-n-1}- t^{n+1}\over (n+1)(t^{-1}-t)}
= {t^{-n} + t^{-n+2} + \cdots + t^n\over n+1}.
\label{potlim}
\end{equation}
For $s\in \mminv$,
\begin{equation}
|\phi(s)| \le 1.
\label{potbound}
\end{equation}
For $s\in 
\C \kern 1pt \backslash \{ (-\infty,0\kern .5pt] \cup \sigma \}$,
\begin{equation}
|\phi(n\kern .3pt s)|^{1/n} \to e^{u(s)} \quad \hbox{as } n\to\infty,
\label{potlim2}
\end{equation}
where
\begin{equation}
u(s) = |\log |t|\kern 2pt | > 0,
\label{potfor}
\end{equation}
with uniform convergence on compact subsets of\/
$\C \kern 1pt \backslash \{ (-\infty,0\kern .5pt] \cup \sigma \}$. 
\end{lemma}

\medskip

\begin{proof}
The formulas (\ref{potbound})--(\ref{potfor}) follow
readily from (\ref{potlim}), so (\ref{potlim}) is what must be established.
First we note that each complex number $s \ne n\kern .7pt\sigma$
corresponds under (\ref{potent}) to two numbers $t$ and $t^{-1}$, and
the right-hand side of (\ref{potlim}) gives the same value for
both, so (\ref{potlim}) does indeed define a function of $s$.
This is a meromorphic function with $n$ zeros at the points
$\alpha_j$ and a pole of order $n$, it can be verified, at
$s=n\kern .7pt \sigma,$ corresponding to
$t=0$.  But at $s=\infty$, corresponding
to $t=1$, this function is analytic.   Thus it is
a degree $n$ rational function with the same
zeros and poles as $\phi(s)$ as defined by (\ref{e7}), and
since it takes the value $1$ at $s=\infty$, it must
be the same function.
\end{proof}

\begin{figure}
\begin{center}
\vspace{10pt}
\indent\hbox{\kern .5in}\includegraphics[scale=.92]{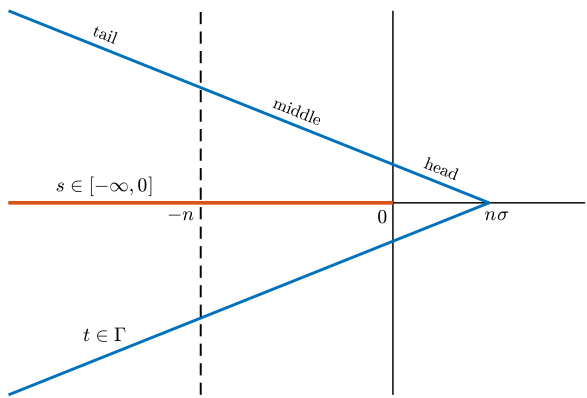}
\vspace{2pt}
\end{center}
\caption{\label{figV}Hankel contour $\Gamma$ for
the proof of Theorem~$\ref{thm2}$ with
$E = \mminv$, with all poles at $s=n\kern .5pt \sigma$.
In the tail,
$F(t)$ in $(\ref{herm})$ is exponentially small.  In the middle,
$F(t)$ is at most algebraically large and
$1/\phi(t)$ is exponentially small. In the head, $F(t)$ is
exponentially large, but its exponential growth rate is limited by $\sigma$, 
so if $\sigma$ is small enough, $F(t)/\phi(t)$ is still exponentially
small.}
\end{figure}

\medskip
{\em Proof of Theorem~$\ref{thm2}$.}
As mentioned above, we take $E = \mminv$ without loss of
generality, and we choose $\Gamma$ for Lemma~\ref{hif}
as a V-shape passing
through $n\kern .5pt \sigma$ at a fixed acute angle,
as sketched in Figure~\ref{figV}.
(Any acute angle will suffice.)
We divide $\Gamma$ into a ``head'' in the
right half-plane, a ``tail'' with $\Re t \le  -n$, and a ``middle''
with $-n \le \Re t \le 0$,
and show that their contributions $\Ihead$, $\Itail$, and $\Imiddle$ to 
(\ref{herm}) are each exponentially small.
(A quantity depending on $n$ is {\em exponentially
small\/} if it is $O(\exp(-C n))$ as $n\to\infty$ for some $C>0$
and {\em exponentially large} if its reciprocal is exponentially small.)
Note that the denominator $t-s$ of (\ref{herm}) is bounded below in absolute
value, $|\phi(s)|$ is bounded by $1$ by (\ref{potbound}),
and $2\pi$ is just a constant.
Thus it is enough to show that each of the three integrals
\begin{equation}
 \int {|F(t)|\over |\phi(t)|} \kern 1pt dt
\label{toshow}
\end{equation}
is exponentially small.
In the tail, $F(t)$ is exponentially
small and decreases exponentially as $\Re t \to - \infty$
by (\ref{h3}), whereas $1/|\phi(t)|$ is bounded by Lemma~\ref{philemma};
thus $\Itail$ is exponentially small.
In the middle, $F(t)$ grows at most algebraically
with $n$ by (\ref{h3}), whereas $1/|\phi(t)|$ is exponentially
small by Lemma~\ref{philemma}; thus
$\Imiddle$ is exponentially small.
In the head,
$F(t)$ is exponentially large, but by (\ref{h4}), at most
$O(\exp(A\kern .7pt \sigma n))$.  On the other hand, by
Lemma~\ref{philemma} again, $1/|\phi(t)|$ is
exponentially small.  In particular, $|\phi(t)| > 
\exp(d\kern .3pt n)$ uniformly 
in this region for some $d>0$ that is independent of $\sigma$.
It follows that if $\sigma < d/A$, then $\Ihead$ too is exponentially small. 
\endproof

\begin{figure}
\begin{center}
\vspace{5pt}
\indent\hbox{\kern .5in}\includegraphics[scale=.92]{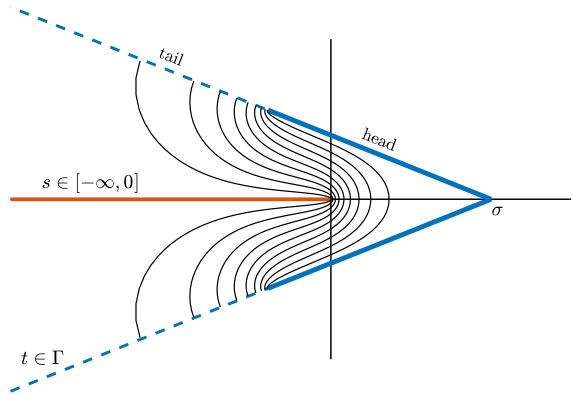}
\vspace{0pt}
\end{center}
\caption{\label{figV2}Hankel contour $\Gamma$ for
the proof of Theorem~$\ref{thmbounded}$ with $E=\inv$,
together with level curves $u = 0.1, \dots, 0.9$ for
the associated potential theory problem obtained
from the conformal map $(\ref{conf1})$--$(\ref{conf2})$.
Poles are distributed according to $(\ref{polechoice})$
along segments of length $\sigma L = O(n)$ on each side of\/ $\Gamma$.
In the tail,
$F(t)$ in $(\ref{herm})$ is exponentially small.  In the head,
$\phi(s)/\phi(t)$ is near-exponentially small.}
\end{figure}

\medskip
We now turn to the second theorem about approximation near $z=0$,
assuming analyticity only in a bounded domain.  As before, our
prototypical example for discussion without loss of generality
is $E = \inv$.  For the argument we again use a V-shaped contour
$\Gamma$ for the Hermite integral, as shown in Figure~\ref{figV2},
but now, $\Gamma$ must be fixed rather than growing with~$n$.
This necessitates a change in the choice of poles.  If the poles
are all at the apex, then the same argument as before ensures
exponentially small contributions to the error from $\Ihead$
and $\Imiddle$, but no longer so small (just root-exponential)
from $\Itail$.  Instead it is necessary to distribute poles along
$\Gamma$, and a successful strategy is to distribute them along
segments extending a distance $\sigma L = \sigma \rho \kern 1pt
n$ from the apex on the two sides of $\Gamma$, for any constant
$\rho> 0$.

Specifically, one way to distribute poles is by solving a
potential theory problem in the infinite $V$ domain.  As suggested
by the contour lines $u = 0.1,0.2,\dots, 0.9$
in Figure~\ref{figV2}, let $u$
be the unique harmonic function in the infinite slit wedge taking values
$u=0$ on $\mminv$ and $u=1$ on the segments of length $\sigma L$ just
mentioned, with a homogeneous Neumann condition on the remainder
of the boundary.  This function $u$ can be calculated by 
means of a conformal map of the lower half of the 
slit wedge onto a rectangle $-K < \Re w < K$, $0 < \Im w < K'$ in
the complex $w$-plane, with the Dirichlet boundary
components mapping to the vertical sides and the Neumann
components to the top and bottom; the ratio $K'/2K$ is known as the
{\em conformal module} of the rectangle.
If the interior angle of
the wedge is $\pi/\mu$ for some $\mu > 1$,
the map has the explicit representation
\begin{equation}
s(w) = \sigma - \sigma \left( {b + y\over 1 + 2b - y}\right)^\mu, \quad
y = \hbox{sn} (w \kern 1pt | \kern 1pt M^{-2}),
\label{conf1}
\end{equation}
where $\hbox{sn}$ is the Jacobi elliptic sine function
and
\begin{equation}
M = 1+2b, \quad  b = L^\mu + \sqrt{L^\mu + L^{2\mu}}.
\label{conf2}
\end{equation}
The numbers $K$ and $K'$ are the complete elliptic
integrals of the first kind with parameters $M^{-2}$ and
$1-M^{-2}$, respectively~\cite{nist}.
Our Fej\'er--Walsh choice of poles and interpolation points
for the theorem will be
\begin{equation}
s_k = s(w_k)\quad w_k = K + i (k-\textstyle{1\over 2})K'/n,
\quad 1 \le k \le n
\label{polechoice}
\end{equation}
and
\begin{equation}
\alpha_k = s(\tilde w_k)\quad \tilde w_k = -K + i (k-\textstyle{1\over 2})K'/n,
\quad 1 \le k \le n.
\label{interpchoice}
\end{equation}

\medskip

\begin{theorem}
\label{thmbounded}
Let\/ $E$ be a simply-connected compact set in 
$\C\kern 1pt \backslash\kern 1pt (-\infty,0\kern .5pt )$; more
generally, $E$ might be a multilevel surface wrapping around\/ $0$
a finite number of times.
Let\/ $f$ be an analytic function in $E$ with 
a branch point singularity at $z=0$ in the
sense of Definition~$\ref{def}$.
Then there exist functions
\begin{equation}
g(z) = c_0^{} + \sum_{k = 1}^n {c_k \over \log(z) - s_k}
\label{ms1}
\end{equation}
satisfying
\begin{equation}
\|\kern .7pt f - g\| \le \exp(- C n / \log n ) \quad \forall n
\label{e3b}
\end{equation}
for some\/ $C>0$, where $\|\cdot\|$ is the supremum norm on $E$.
If $E = \inv$, a suitable choice of
$\{s_k\}$ is $(\ref{polechoice})$ for any\/ $\mu$ with $\tan(\pi/2 \mu)
< \tau$ and sufficiently small\/ $\sigma$.
\end{theorem}
\medskip

\begin{proof}
As with Theorem~\ref{thm1}, the proof is carried out via the
Hermite integral formula for the equivalent problem of rational
approximation of a function $F(s)$ on the domain $\log(E)$ in
the $s = \log z$ variable; again we first remove $f(0)$ from
the problem by setting $c_0^{} = f(0)$.  If $E$ is a compact
set in $\C\kern 1pt \backslash\kern 1pt (-\infty,0\kern .5pt )$,
or a multilevel surface wrapping around\/ $0$ a finite number of
times, then $\log (E)$ is a closed subset of $\C \cup \{\infty\}$
with upper-bounded\ $\Re E$ and\/ $|\Im E|$.  We take
$E = \inv$ and $\log(E) = \minv$ for simplicity, as sketched
in Figure~\ref{figV2}; the more general case is treated by
an adjustment of the contour $\Gamma$ and the conformal map.
By the assumption (\ref{taucond}), $F$ can be extended to a
single-valued analytic function throughout a V-shaped region in
the $s$-plane, as in the figure, for any inner angle $\pi/\mu$ with
$\tan(\pi/2\mu) < \tau$ and sufficiently small $\sigma$.  As in the
proof of Theorem~\ref{thm1}, we now estimate the Hermite integral
(\ref{herm}) over the ``head'' region of $\Gamma$, consisting of
the two segments of length $\sigma L = \sigma \rho \kern 1pt n$
touching the apex, and the remaining ``tail.''  The contribution
$\Itail$ is exponentially small because of (\ref{bpdef}),
which in the $s$ variable becomes (\ref{h3}).  The contribution
from $\Ihead$ is near-exponentially small, which we establish as
follows.  The factor $F(t)/(t-s)$ in (\ref{herm}) is bounded as a
consequence of (\ref{h3}) and the contour $\Gamma$ being fixed,
so we need to show that $\phi(s)/\phi(t)$ is near-exponentially
small.  From potential theory related to Theorem~19 of chapter 9
of~\cite{walsh}, it is known that this selection of $\{s_k\}$ and
$\{\alpha_k\}$ as Fej\'er--Walsh points ensures $|\phi(s)/\phi(t)|
\approx \exp(-n(K/K'))$ as $n\to\infty$.  The result follows from
the estimate $K'/K \sim (4\mu / \pi) \log n$ as $n\to\infty$,
which can be derived from~\cite[eq.~(19.9.5)]{nist}.  \end{proof}

\section{\label{sec5}Convergence theorem: multiple singularities}

Theorems~\ref{thm1} and~\ref{thmbounded} are local assertions,
establishing exponential and near-exponential resolution
of isolated singularities.  It remains to show that global
approximations can be constructed by adding together these local
pieces.  We do this in the style of Theorem~\ref{thmbounded},
requiring analyticity just in a neighborhood of $E$.  The
argument is adapted from the discussion around Theorem~2.3
in~\cite{lightning} and illustrated schematically in
Figure~\ref{fig6}.

\medskip

\begin{theorem}
\label{thmmult}
Let\/ $E$ be a simply-connected compact set in 
$\C$ and let $f$ be an analytic function in $E$ with
branch point singularities in the sense of Definition~$\ref{def}$ at
boundary points $z_1,\dots, z_m$; more generally, 
$E$ might be a multilevel surface wrapping around each
branch point a finite number of times.
Then there exist functions
\begin{equation}
g(z) = c_0^{} + \sum_{j = 1}^m \kern 2pt \sum_{k=1}^{n_j}\kern 1pt
{c_{jk} \over \log(z-z_j) - s_{jk}}
\label{formsmult}
\end{equation}
satisfying
\begin{equation}
\|\kern .7pt f - g\| \le \exp(- C n / \log n ) \quad \forall n
\label{e3mult}
\end{equation}
for some\/ $C>0$, where $\|\cdot\|$ is the supremum norm on $E$.
\end{theorem}
\medskip

\begin{proof}
The function $f$ can be represented on $E$ as a Cauchy
integral
\begin{equation}
f(z) = {1\over 2\pi i} \int_\Gamma {f(t)\over t-z} \kern 1pt dt, 
\label{cauch}
\end{equation}
where $\Gamma$ is any fixed contour that lies outside $E$ and within
the region of analyticity of $f$, except that $\Gamma$ touches
$E$ at each of the points $\{z_k\}$.
We split up $\Gamma$ into $m$ pieces $\Gamma_k$, with $\Gamma_k$ touching
just $z_k$, giving
\begin{equation}
f(z) = f_1(z) + \cdots + f_m(z),
\end{equation}
where $f_k$ is the Cauchy integral evaluated over the arc
$\Gamma_k$,
\begin{equation}
f_k(z) = {1\over 2\pi i} \int_{\Gamma_k^{}} {f(t)\over t-z} \kern 1pt dt.
\label{cauchk}
\end{equation}
For the mathematics of Cauchy
integrals over open arcs, see Chapter 14 of~\cite{acca}.

\begin{figure}
\begin{center}
\vspace{9pt}
\indent\hbox{\kern .5in}\includegraphics[scale=1]{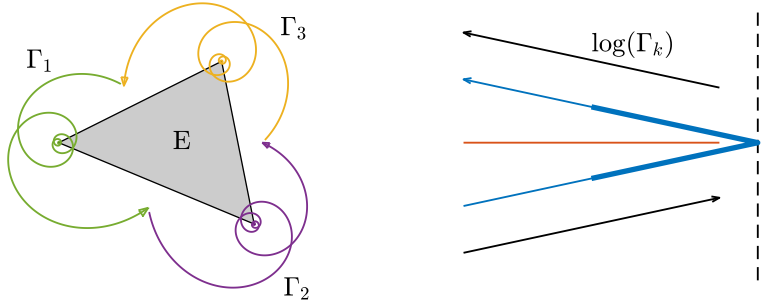}~~~~~~~
\vspace{0pt}
\end{center}
\caption{\label{fig6}Illustrations for the proof
of Theorem~$\ref{thmmult}$.  A function $f$ on $E$ is decomposed
into $m$ pieces, with $f_k$ defined by a Cauchy integral of\/ $f$ over an arc 
$\Gamma_k$ touching the branch point $z_k$;
$\Gamma_k$ is a pair of logarithmic spirals on a Riemann
surface to which
$f$ is extended by analytic continution.
The right image shows the configuration near a
branch point after a change of variables
$s = \log(e^{i\alpha}(z-z_k))$.  The arc defining $f_k$
becomes a disjoint pair of rays $\log(\Gamma_k)$
in the $s$-plane, and poles
$\{s_k\}$ of the rational approximation
will be placed along segments of length $O(n)$
on a V-shaped contour lying
between $\log(\Gamma_k)$ and $\mminv$.}
\end{figure}

We are done if we can show that each function $f_k$ can
be approximated on $E$ with near-exponential accuracy by a
function of the form (\ref{ms1}), with $\log(z)$ replaced by
$\log(z-z_k)$ (that is, a continuous branch of $\log(z-z_k)$
in $E$).  This will be ensured by Theorem~\ref{thmbounded} if
$f_k$ satisfies the conditions (\ref{taucond}) and (\ref{bpdef})
in our definition of a branch point singularity.  At this stage,
the choice of $\Gamma$ becomes crucial.  A choice of $\Gamma$
as a closed, non-intersecting contour in the $z$-plane will
not do; the corresponding functions $f_k$ will not be defined
along curves winding around $z_k$ as required by (\ref{taucond}).
Instead, $\Gamma$ must be taken to be a spiral on a Riemann
surface for $f$ whose projection on the $z$-plane is a self-intersecting contour
with logarithmic spiral behavior at each $z_k$, as sketched in
Figure~\ref{fig6}.  The integrals defining each $f_k$ now make
sense since each $z_k$ is a branch point around which $f$ can be
analytically continued.  Condition (\ref{bpdef}) holds
since it holds for $f$ by assumption and $f_k$ differs from $f$ by
a function that is analytic near $z_k$, namely the sum of
the contributions $f_j$ with $j\ne k$.
\end{proof}

\section{\label{secsing}Behavior near singularities}
\begin{figure}
\begin{center}
\vspace{8pt}
\includegraphics[scale=.95]{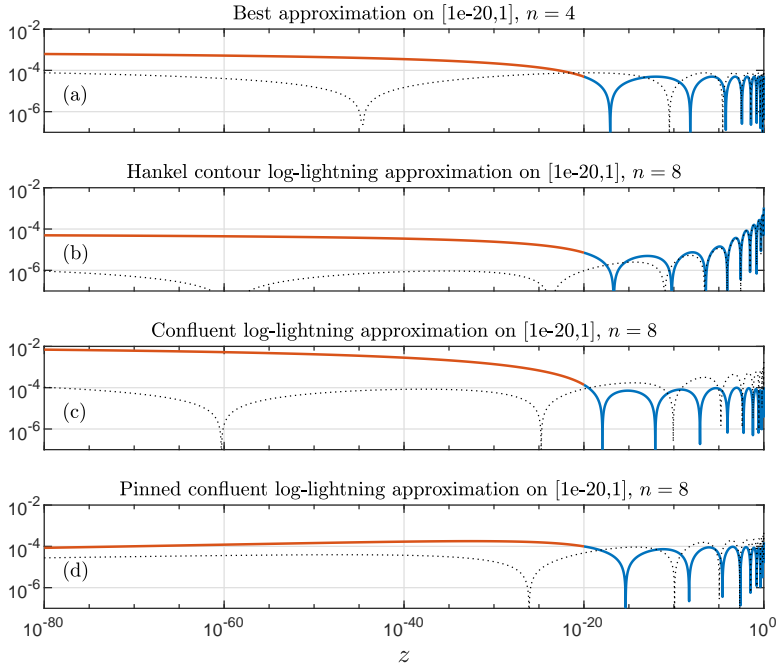}
\vspace{-3pt}
\end{center}
\caption{\label{fig7}Four approximations of $\sqrt z$ computed on
$\smallint$ and then plotted over
$\bigint$ (the absolute value of the error). 
The fine dots show corresponding results
for an approximation computed over $\bigint$.  See the text
for discussion.}
\end{figure}

As mentioned in the penultimate paragraph of section~\ref{sec2}
and illustrated in Figure~\ref{fig3}, reciprocal-log approximations
have surprising properties near the singularities.  Whereas error
curves for rational approximations vary over a scale exponentially
close to the singularity, for reciprocal-log approximations this
becomes doubly-exponential.  Fortunately, so far as we are aware,
these effects need not cause difficulties in using
these approximations, and in particular, it appears that they
do not necessitate the use of extended-precision arithmetic.
We shall now explain our understanding of these matters with
reference to Figure~\ref{fig7}, which presents four
approximations of $\sqrt z$ for $z\in \inv$.  In each case
the approximation is computed over $\smallint$ and then the
absolute value of the error is plotted over $\bigint$.
For comparison, fine dots show the corresponding results for
computations over $\bigint$.

Image (a) shows the best (minimax) reciprocal-log approximation
with $n=4$.  We see here that the best approximation over
$\smallint$ is by no means optimal over $\bigint$.  The error
is 10 times larger there (6.2e-4 vs.\ 5.0e-5), and this ratio
will grow rapidly with increasing $n$.  One might think this
implies that reciprocal-log approximations must be ineffective
at approximating all the way up to the singularity, but the small
dots in the figure, corresponding to the best approximation over
$\bigint$, contradict this expectation, showing an error of
just 7.6e-5.  As can be seen in Figure~\ref{fig3},
the error in best approximation over all of $\inv$
is also not much larger, just 8.7e-5.

Image (b) shows the error for reciprocal-log approximation
(\ref{partfracr}) with poles positioned on a Hankel contour
according to (\ref{hankelpts}), now with $n=8$.  Again the error
in the left-hand part of $\bigint$ is larger than it needs to be,
but this may be less important for these approximations since
the error is dominated by larger values of $z$.

Image (c) shows a very different situation for reciprocal-log
approximation (\ref{partfracrc}) with confluent poles at $s_0
= n/2$.  Now the approximation over $\smallint$ gives an error
over $\bigint$ that is orders of magnitude too large.

It is at this point that we wish to pause and ask, what might be
the practical implications of such behavior?  It certainly seems
disturbing if an approximation is much less accurate on $\bigint$
than on $\smallint$.  There are three rather disparate observations
to be made, which combine to an encouraging picture.

First we note that in IEEE floating point arithmetic, whereas
one can represent numbers near $0$ smaller than $10^{-300}$,
the resolution near other points $z_k$ is only on the order of
16 digits.  Thus even if we wanted to compute approximations by
least-squares fitting on a grid in $z_k + (0,10^{-20}\kern .8pt
]$, we could not do so.  This might seem worrying, if one takes
Figure~\ref{fig7} to suggest that such domains will be required.

The second observation is that as a purely practical matter,
they should not be required.  On a planar domain $E$ such as
that of Figure~\ref{figphase}, for example, if a function $f$
is approximated to a satisfactory precision everywhere except
at distances $< 10^{-20}$ from the vertices $\{z_k\}$, what
difference does it make?  What is wrong with an approximation
that is mathematically inaccurate in theory, but only at points
where it will never be evaluated?

And yet the use of such approximations would be troubling.
This brings us to the third observation, which is the basis of
image (d) of Figure~\ref{fig7}: it may be possible to make an
approximation accurate on $z_k + (0,10^{-20}\kern .8pt]$ without
doing any arithmetic there.  The trick is to forget least-squares
gridding in $z_k + (0,10^{-20}\kern .8pt]$ and instead impose
additional conditions at $z_k$, forcing the singular part of
the approximation to decay there at a rate $O(1/(\log|z-z_k|)^J)$
for some power $J>1$.  This is natural since in floating point
arithmetic, all of $z_k + (0,10^{-20}\kern .8pt]$ rounds to
$z_k$ anyway.  Image (d) corresponds to an approximation of
this kind with $J=2$: we approximate in the $n$-dimensional space spanned
by $1/(\log(z)-n/2)^j$ with $2\le j \le n+2$.  The work and the
number of degrees of freedom are unchanged, and now least-squares
fitting on $\smallint$ gives an approximation that is accurate
on all of $\inv$.  We call this a {\em pinned\/} approximation,
since certain terms of the approximation are constrained to be zero
at the singularity.  In this example $J=2$ is effective,
but for larger $n$ a choice closer to $J=n/2$ will be called for.
To estimate this number, one can use analysis such as that of
(\ref{interppts}) or (\ref{interpchoice}) to monitor how many
interpolation points might be expected to fall closer to $z_k$
than about machine precision.  For the problem of image (d)
with $n=8$, (\ref{interppts}) gives distances $\exp(\alpha_j)
\approx 0.8, 0.6, 0.3, 0.06, 0.003, 10^{-5}, 10^{-13}, 10^{-56}$,
and our pinned approximation might be thought of as effectively
replacing the last two of these numbers by zero.  There are many
mathematical questions here and it would be interesting to
investigate them.

\section{\label{sec7}Log-lightning PDE solvers}

The motivation for this work has been the development of
numerical methods for the solution of partial differential
equations on planar domains with corners, beginning with the
Laplace equation with Dirichlet boundary conditions.  For these
problems we do not know an analytic function on $E$, just its
real part on the boundary.  This is no difficulty for linear
least-squares fitting, however, where all that matters is having
an efficient approximation space.  The idea of solving Laplace
problems in this fashion was presented in~\cite{lightning},
where root-exponential convergence of approximations based on
rational functions is established theoretically and experimentally.
Here we move to reciprocal-log approximations, and the theorems
of sections~\ref{sec4} and~\ref{sec5} assert that the convergence
should improve to exponential or near-exponential.

\begin{figure}
\begin{center}
\vspace{6pt}
\includegraphics[scale=1.10]{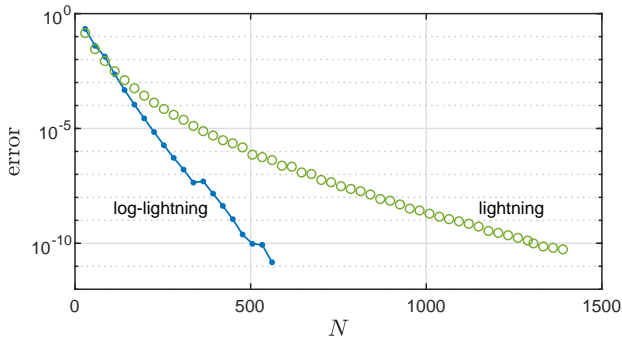}
\vspace{-3pt}
\end{center}
\caption{\label{fig8}Convergence of lightning and log-lightning solutions
to the NA Digest problem of the Laplace equation on an L-shaped
region~{\rm\cite{digest}}.  The difference between root-exponential
and exponential convergence is evident.}
\end{figure}

Figure~\ref{fig8} presents a pair of curves confirming this
expectation.  The solutions computed here are for the problem
posed on NA Digest in December 2019, defined by the L-shaped region
with vertices $0, 2, 2+i, 1+i, 1+2i, 2i$ and boundary data $u(z)
= \Re (z^2)$~\cite{digest}.  As usual, $N$ denotes the total
number of degrees of freedom, which in this plot is a number of
the form $14 n+1$ with $n = 2,4,6,\dots .$  For each $n$, the
approximating function consists of the real part of a polynomial
of degree $n$ together with $n$ singularities $\{s_k\}$ at each
of the six corners as defined by (\ref{hankelpts}), though with
a factor $n/3$ instead of $n/4$ in front.  This gives $7(n+1)$
complex degrees of freedom, hence $14n + 2$ real ones, and the
number reduces to $14n+1$ since the imaginary part of the constant
term of the polynomial does not affect the fit.

This experiment is a long way from the software that has been
developed for lightning approximation~\cite{laplace}.  The boundary
has been resolved simply by 500 exponentially clustered points on
each side, whereas adaptive software can make the gridding both
more efficient and more careful.  More importantly, each corner
has been allocated the same number $n$ of singularities, whereas
adaptive software will put more singularities at vertices with
stronger singularities.  So although the log-lightning convergence
curve of Figure~\ref{fig8} is promising, the reliable fast way
to solve Laplace problems at this date is still by means of the
{\tt laplace} code at~\cite{laplace}.  This code also addresses
a number of variations such as Neumann boundary conditions,
discontinuous data, and piecewise smooth curved boundaries.

The number of degrees of freedom $N$ is not a direct measure of
computer time.  For reasons of linear algebra, the work actually
increases in proportion to $N^2$ or $N^3$, depending on whether
the grid is fixed or refined with $N$, which tends to make the
difference between the lightning and log-lightning computations
greater than it may appear in Figure~\ref{fig8}.  A consideration
pushing in the other direction is that computer evaluation of
a complex logarithm $\log z$ is typically slower than that of a
reciprocal $1/z$ by a factor such as $3$ or $4$.

It was shown in \cite{helm} that the original lightning Laplace
solver can be generalized to a lightning Helmholtz solver by
replacing poles by Hankel functions.  Presumably there is a
Helmholtz generalization of the log-lightning method too, but we
have not investigated this.

\section{\label{sec8}Conclusion}
Exponential or near-exponential convergence of approximations to
branch point singularities is a new phenomenon.  Reciprocal-log
approximations are interesting in theory and may have consequences
in practice if good software can be developed.

Philosophical questions are attached to this kind of approximation.
It is an old idea that polynomials and rational functions have
a special status since ``the only operations that can really
be carried out numerically are the four elementary operations of
addition, subtraction, multiplication and division.''\footnote{This
quote comes from the PhD thesis of Hilbert's student Kirchberger
in 1902~\cite[pp.~196--197]{atap}.}  This point of view has
contributed to the dominance of polynomials and rational functions
in approximation theory, with other classes of approximating
functions seeming less fundamental.  We are not ourselves immune
to the impression that there is something contrived about
approximations of the forms (\ref{form}) and (\ref{formsm}).
But the philosophical distinction that may have seemed sharp
before the arrival of computers seems less sharp today.

A fascinating aspect of reciprocal-log approximations is their
exploitation of the behavior of a function on a Riemann surface.
As one of the smallest consequences of this feature, these
approximations have no difficulty at all in treating domains with
slits, and the prospects for further adventures with multivalued
functions seem enticing.

\section*{Acknowledgments}
The research for this article was carried out in regular
communication with Andr\'e Weideman, whose many suggestions
we are happy to acknowledge.  Another key contribution was
Lemma~\ref{philemma}, from Peter Baddoo.  In addition we are
grateful for helpful comments from Marco Fasondini, Abi Gopal,
Alex Townsend, and Elias Wegert.

\indent~~\kern 1in

\end{document}